\documentclass[12 pt, draftcls, onecolumn]{IEEEtran}
\usepackage{cite}

\ifCLASSINFOpdf
\else
\fi
\usepackage{amsmath}
\usepackage{algorithmic}

\usepackage{array}
\usepackage{fixltx2e}

\usepackage{stfloats}

\usepackage{algorithm}
\usepackage{algorithmic}
\usepackage{amsfonts}
\usepackage{graphicx}
\newtheorem{lemma}{Lemma}
\newtheorem{remark}{Remark}
\newtheorem{example}{Example}
\newtheorem{assumption}{Assumption}
\newtheorem{definition}{Definition}
\newtheorem{theorem}{Theorem}
\newcommand{\R}{\mathbb{R}}
\newcommand{\rr}{\mathcal{R}}
\newcommand{\cc}{\mathcal{C}}

\begin{document}
%
\title{Distributed Computation for Solving the Sylvester Equation Based on Optimization}
%
%
%

\author{Wen~Deng,~
        Xianlin~Zeng,~
        and~Yiguang~Hong,~\IEEEmembership{Fellow,~IEEE}
\thanks{W. Deng and Y. Hong are with the School of Mathematical Sciences, University of Chinese Academy of Sciences; Key Laboratory of Systems and Control, Academy of Mathematics and Systems Science, Chinese Academy of Sciences, Beijing, 100190, China (dengwena@amss.ac.cn, yghong@iss.ac.cn).}
\thanks{X. Zeng is with Key Laboratory of Intelligent Control and Decision of Complex Systems, School of Automation,
	Beijing Institute of Technology, Beijing, 100081, China (xianlin.zeng@bit.edu.cn).}}

%
%

\markboth{Journal of \LaTeX\ Class Files
}%
{Deng \MakeLowercase{\textit{et al.}}: Distributed Computation for the Sylvester Equation}
%



\maketitle

\begin{abstract}
This paper solves the Sylvester equation in the form of $ AX+XB=C $ in a distributed way, and proposes three distributed continuous-time algorithms for three cases.  We start with the basic algorithm for solving a least squares solution of the equation, and then give a simplified algorithm for the case when there is an exact solution to the equation, followed by an algorithm with regularization case.  Based on local information and appropriate communication among neighbor agents, we solve the distributed computation problem of the Sylvester equation from the optimization viewpoint, and we  prove the convergence of proposed algorithms to an optimal solution in  three different cases, with help of the convex optimization and semi-stability.
\end{abstract}


%
\IEEEpeerreviewmaketitle

\section{Introduction}
%
%
%
%
\IEEEPARstart{U}{nder}   the influence of big data and large-scale systems, distributed optimization and computation have attracted more and more research attention.
Both discrete-time algorithms \cite{nedic2009distributed,wang2017distributed,yuan2016regularized} and continuous-time algorithms \cite{gharesifard2014distributed,jafarzadeh2017optimizing,deng2019distributed} have been given for various distributed optimization problems. The basic idea is that many interconnected agents in a network, having local information separately, cooperate with their neighbors to exchange information and achieve global goals eventually.

With the rapid development of distributed optimization algorithms, the idea of using distributed methods to solve matrix equations has attracted much interest. The Sylvester equation is an important class of matrix equations, which has wide application in control theory, systems theory and many other fields \cite{datta2004numerical,bhatia1997how,simoncini2016computational,corless2003linear}.
For instance, the Sylvester equation plays a significant role in computing invariant subspaces \cite{demmel1987three},  achieving pole assignment \cite{lv2010parametric} and model reduction \cite{sorensen2002sylvester}.
There have been many centralized algorithms for solving matrix equations, such as Schur decomposition methods, Krylov-subspace methods, and iterative methods \cite{simoncini2016computational,hu1992krylov,benner2009adi,ding2006iterative}.
However, those centralized algorithms for solving the Sylvester equation $ AX+XB=C $ mainly need to deal with the whole two coefficient matrices $ A $ and $ B, $ which could not be applied directly to many distributed scenarios, including the one we consider in this paper. Besides, the parallel distributed computation for the Sylvester equation \cite{granat2003parallel} could not work only for such local information either, because it needs to transform $ A $ and $ B $ to real Schur form at first.

Recently, there have been several works about solving the linear algebraic equation in the form of
$ Ax=b $ by various distributed methods  \cite{shi2017networkTAC,liu2019arrow,liu2017asynchronous,wang2017double,cao2017continuous}. With the help of a network structure, every node only gets access to the local information, for example, a row \cite{shi2017networkTAC,liu2019arrow,liu2017asynchronous} or a column \cite{cao2017continuous} of the matrix $ A $. In \cite{wang2017double}, there is a double-layered framework for all nodes and it enhances the flexibility of the access to information. There is not much work about solving matrix equations in a distributed way yet. Though we could transform a matrix equation into the form of a linear algebraic equation sometimes, its partition structure has less flexibility to a certain extent. Due to matrix multiplication rules, there are some differences between matrix equations and linear algebraic equations. A class of matrix equations formed as $ AXB=F $ was first discussed in \cite{zeng2018distributed}, with different distributed algorithms according to different partition structures.
However, there are no results on the distributed computation of the Sylvester equation, which is more complicated than $ AXB=F $, to our knowledge.

The objective of this paper is to solve the Sylvester equation in a specific distributed formulation.  Notice that these methods in $ Ax=b $ \cite{shi2017networkTAC,liu2019arrow,liu2017asynchronous} and $ AXB=F $ \cite{zeng2018distributed} cannot be applied directly to the Sylvester equation $ AX+XB=C $ considered in this paper, because properties of matrix multiplication and the given partition structure make it different from $ Ax=b $ and $ AXB=F $.
In the  distributed structure we consider, each agent only needs to know partial rows or columns of information matrices $ A, B, C $ and communicate with its  neighbors to exchange information.
Therefore, appropriate conversion methods are introduced  to deal with the given partition structure of $ AX+XB=C $ in our problem formulation.
Main contributions of this paper are summarized as follows:
\begin{itemize}
	\item We consider the distributed computation problem for the Sylvester equation with a special distributed structure in two main cases: the least squares solution case and the regularization case.  To design the distributed algorithm, in which each agent exactly knows the corresponding rows of $ A $ and columns of $ B $,  we introduce two new variables to deal with the inconsistency caused by different partitions of rows or columns, and then reformulate the problem as some distributed optimization problems by constructing an appropriate equivalent transformation.
	\item We propose distributed algorithms to solve the Sylvester equation in different cases based on saddle-point dynamics, and moreover, we propose a simplified algorithm for the case when there exist exact solutions.   It is worthwhile to mention that we employ the derivative feedback technique in the (nonsmooth) regularization case for the optimization problem with the penalty term for variable $ X $, to solve the problem.
	\item
		We provide exponential convergence analysis of proposed algorithms for the least squares (exact) solution  case and the regularization case, respectively, under mild conditions. Since solutions of the Sylvester equation may be not unique and the optimal solution obtained by the algorithm trajectory may depend on initial conditions, we employ the semi-stability theory of  dynamical systems, which tackles convergence properties of dynamical systems having a continuum of equilibria.
\end{itemize}

The rest paper is organized as follows. Section \ref{sec:math pre} introduces  relevant mathematical preliminaries, while Section
\ref{sec:problem formulation} formulates  distributed computation problems of the Sylvester equation and  reformulated distributed optimization problems.
Then Sections \ref{sec:exact+LS} and \ref{sec:regularization} present distributed algorithms for the least squares solution case, the exact solution case, and the regularization case with corresponding convergence analysis.
Following that, Section \ref{sec:simulation} provides a short discussion and three numerical examples, and Section \ref{sec:conclusion} concludes this paper briefly.

\section{Preliminaries}
\label{sec:math pre}
In this section, we give related notations,  some basic concepts, and  some lemmas  for the analysis.

 Let $ \mathbb{R}^{m\times n} $ denote the set of $ m\times n $ real matrices, $ \mathbf{1}_n~( \mathrm{or}~\mathbf{0}_n)$ and $ I_n $ denote a vector in $ \mathbb{R}^n $ with all elements of 1 (or 0) and the $ n\times n $ identity matrix, respectively. Let $ M^T, \mathrm{Im}(M) $, $ \mathrm{ker}(M) $ and $ \mathrm{spec}(M) $ denote the transpose, the image, the kernel and the set of all eigenvalues of the matrix $ M $, respectively. Let $ \mathrm{vec}(X)\in\mathbb{R}^{mn} $ denote the vector that is a stack of all columns in matrix $ X\in\mathbb{R}^{m\times n} $  and col $ \{M_1, \ldots, M_n\}=[M_1^T, \ldots, M_n^T]^T $ for every $ M_i\in\mathbb{R}^{m_i\times n} $. Let $ \otimes $ denote the Kronecker product. Let $ ||\cdot||_F $ denote the Frobenius norm of a matrix, $ ||M||_F:=\sqrt{\sum\nolimits_{i,j}M_{i,j}^2} $ and $ \langle \cdot , \cdot\rangle_F $ denote the Frobenius inner product of real matrices:
$ \langle M_1, M_2\rangle_F=\sum\nolimits_{i,j}(M_1)_{i,j}(M_2)_{i,j} $
with $ M_1, M_2\in \mathbb{R}^{m\times n}. $
Define augmented matrices $ [Y_i]_\rr , [Z_i]_\cc\in\mathbb{R}^{m\times r} $ as
\begin{align*}
&[Y_i]_\rr=\begin{bmatrix} 0_{m_1\times r}^T,   \cdots  0_{m_{i-1}\times r}^T, Y_i^T, 0_{m_{i+1}\times r}^T, \cdots , 0_{m_{n}\times r}^T\end{bmatrix}^T,\\& [Z_i]_\mathcal{C}=\begin{bmatrix}
0_{m\times r_1}, \cdots, 0_{m\times r_{i-1}}, Z_i, 0_{m\times r_{i+1}}, \cdots, 0_{m\times r_n}
\end{bmatrix},
\end{align*}
with $  Y_i\in \mathbb{R}^{m_i\times r}, Z_i\in \mathbb{R}^{m\times r_i}. $

\begin{definition}[see \cite{rockafellar1970convex}]
	Let $ f $ be a convex function. A vector $ h $ is called a subgradient of function $ f $ at point $ X_0\in \mathrm{dom}f $ if, for any $ X\in  \mathrm{dom}f $,
	\begin{equation}\label{g1}
	f(X)\ge f(X_0)+\langle h, X-X_0\rangle.
	\end{equation}
	The set of all subgradients $ \partial f(X_0) $ is said to be the subdifferential of function $ f $ at point $ X_0. $
\end{definition}

Consider a time-invariant dynamical system
\begin{equation}\label{ds}
\dot{x}(t)=\phi(x(t)), \quad x(0)=x_0, \quad t\ge 0,
\end{equation}
where $ \phi: \mathbb{R}^d\rightarrow\mathbb{R}^d $ is Lipschitz continuous.
Given $ x(t): [0,\infty)\rightarrow \mathbb{R}^d $ as a trajectory of system \eqref{ds},
the point $ z\in \mathbb{R}^d $ is a limit point of a solution $ x(t) $ if there exists a positive increasing divergent sequence $ \{t_i\}_{i=1}^{\infty}\subset \mathbb{R} $ such that $ z=\lim_{i\rightarrow\infty}x(t_i) $.
A set $ \mathcal{M} $ is positively invariant with respect to \eqref{ds} if, for every $ x_0\in \mathcal{M}, \mathcal{M} $ contains the solution $ x(t) $ of \eqref{ds} for all $ t> 0 $ with $ x(0)=x_0 $.  Then we introduce the following definition from \cite{hui2009semistability} and Lemma \ref{lemma:converge}, which are useful   for convergence analysis of algorithms in the least squares case.

\begin{definition}\label{def:semi-stable}
	A point	$ z $ is semi-stable if $ z $ is Lyapunov stable and there exists an open subset $ \mathcal{D}_0 $ of $ \mathcal{D} $ containing $ z $ such that, for all initial conditions in $ \mathcal{D}_0 $, the solution of \eqref{ds} converge to a Lyapunov stable equilibrium point.   Then system \eqref{ds} is semi-stable with respect to $ \mathcal{D} $ if every solution with initial condition in $ \mathcal{D} $ converges to a Lyapunov stable equilibrium.  Moreover, \eqref{ds} is said to be globally semi-stable if it is semi-stable with respect to $ \mathbb{R}^d $.
\end{definition}

\begin{lemma}[Theorem 3.1 in\cite{hui2009semistability}]\label{lemma:converge}	
	Let $ \mathcal{D} $ be an open positively invariant set with respect to \eqref{ds}, $ V: \mathcal{D}\rightarrow\mathbb{R} $ be a continuously differentiable function, and $ x(t) $ be a solution of \eqref{ds} with $ x(0)\in \mathcal{D} $, contained in a compact subset of $ \mathcal{D} $.
	Assume $ \frac{d}{dt}V(x(t))\le 0, $ for all $ x\in \mathcal{D} $ and define $ \mathcal{Z}=\{x\in \mathcal{D}: \frac{d}{dt}V(x)=0\}. $ If every point in the largest invariant subset $ \mathcal{M} $ of $ \bar{\mathcal{Z}}\cap \mathcal{D} $ is Lyapunov stable equilibrium, where $ \bar{\mathcal{Z}} $ is the closure of $ \mathcal{Z} $, then the system \eqref{ds} is semi-stable
	with respect to $ \mathcal{D} $.
\end{lemma}

Consider a time-invariant differential inclusion
\begin{equation}\label{eq:different-inclusion}
\dot{x}(t)\in \mathcal{F}(x(t)), \quad x(0)=x_0, \quad t\ge 0,
\end{equation}
where $ \mathcal{F}: \mathbb{R}^d\rightarrow\mathcal{B}(\mathbb{R}^d) $ is upper-semicontinuous, and the set $ \mathcal{F}(x) $ is nonempty, compact and convex for any $x$, which can guarantee the existence of solutions of \eqref{eq:different-inclusion}\cite{cortes2008discontinuous}. A set $ \mathcal{M} $ is weakly (strongly) positively invariant with respect to \eqref{eq:different-inclusion} if, for every $ x_0\in \mathcal{M}, \mathcal{M} $ contains a solution $ x(t) $ (all solutions) of \eqref{eq:different-inclusion} for all $ t> 0 $ with $ x(0)=x_0 $.
 An equilibrium $ x^* $ of \eqref{eq:different-inclusion} satisfies $ \mathbf{0}_d\in \mathcal{F}(x^*) $.  The definition of semi-stable about \eqref{eq:different-inclusion} \cite{hui2009semistability} is similar to Definition \ref{def:semi-stable}. Then we introduce the following Definition \ref{def:lie} and Lemmas \ref{lemma:df-dt} and \ref{lemma:nonsmooth_stable} from \cite{cortes2008discontinuous}, \cite{bacciotti1999stability}, which are useful  for convergence analysis of the algorithm in the regularization case.
\begin{definition}\label{def:lie}
	The set-valued Lie derivative $ \tilde{\mathcal{L}}_{\mathcal{F}}f:\R^d \to \mathcal{B}(\R) $ of a locally Lipschitz function $ f:\R^d\to \R $ with respect to $ \mathcal{F}:\R^d\to \mathcal{B}(\R^d) $ at $ x $ is defined as
	\begin{equation}\label{eq:Lie-derivative}
	\begin{aligned}
	\tilde{\mathcal{L}}_{\mathcal{F}}f(x)=\{&a\in \R: \mathrm{there\ exists}\ y\in \mathcal{F}(x) \mathrm{such\ that}\\& h^Ty=a,\ \forall h\in \partial f(x) \}.
	\end{aligned}	\end{equation}
\end{definition}
\begin{remark}
	If $ f $ is differentiable at $ x, $ then
	$\tilde{\mathcal{L}}_{\mathcal{F}}f(x)=\{\nabla f(x)^Ty:y\in\mathcal{F}(x)\}.$
\end{remark}

\begin{lemma}\label{lemma:df-dt}
	Let $ x(t) $ be a solution of different inclusion \eqref{eq:different-inclusion}, and let $ f:\R^n\to \R $ be a locally Lipschitz continuous and regular function \cite{cortes2008discontinuous}. Then $ \frac{d}{dt}f(x(t)) $ exists and $ \frac{d}{dt}f(x(t))\in\tilde{\mathcal{L}}_{\mathcal{F}}f(x) $ almost everywhere.
\end{lemma}
\begin{lemma}\label{lemma:nonsmooth_stable}
	Suppose that \eqref{eq:different-inclusion} has at least one solution and $ x^* $ is an equilibrium of the \eqref{eq:different-inclusion}. Let $ f $ be locally Lipschitz continuous and regular, and satisfy (\romannumeral1 ) $ f(x)\ge 0, f(x)=0 $ if and only if $ x=x^*; $ (\romannumeral2 ) $ \max\tilde{\mathcal{L}}_{\mathcal{F}}f(x)\le 0. $ Then $ x^* $ is a Lyapunov stable equilibrium of \eqref{eq:different-inclusion}.
\end{lemma}

\begin{lemma}[Theorem 3.1 in\cite{hui2009semistability}]\label{lemma:regular-converge}
	Let $ \mathcal{D} $ be an open strongly positively invariant set with respect to \eqref{eq:different-inclusion}, $ V: \mathcal{D}\rightarrow\mathbb{R} $ be a locally Lipschitz continuous and regular function, and $ x(t) $ be a solution of \eqref{eq:different-inclusion} with $ x(0)\in \mathcal{D} $, contained in a compact subset of $ \mathcal{D} $.
	Assume $ \tilde{\mathcal{L}}_{\mathcal{F}}V(x)\le 0 $ for all $ x\in \mathcal{D} $ and define $ \mathcal{Z}=\{x\in \mathcal{D}: 0\in \tilde{\mathcal{L}}_{\mathcal{F}}V(x)\}. $ If every point in the largest weakly positively invariant subset $ \mathcal{M} $ of $ \bar{\mathcal{Z}}\cap \mathcal{D} $ is a Lyapunov stable equilibrium, where $ \bar{\mathcal{Z}} $ is the closure of $ \mathcal{Z} $, then the system \eqref{eq:different-inclusion} is semi-stable
	with respect to $ \mathcal{D} $.
\end{lemma}

\section{Problem Formulation}
\label{sec:problem formulation}
In this section, we first formulate our problem, and then introduce some transformations which will be helpful in the distributed reformulation.

Consider a continuous-time Sylvester equation of the following form
\begin{equation}\label{se}
AX+XB=C,
\end{equation}	
with $ A\in \mathbb{R}^{m\times m}, B\in \mathbb{R}^{r\times r}, C\in \mathbb{R}^{m\times r} $ and the unknown variable $ X\in \mathbb{R}^{m\times r} $. Here we consider a specific partition for $ n $ parts, Left-Row-Right-Column (LRRC) case, to solve the Sylvester equation: dividing matrix $ A $ by row and matrix $ B $ by column. Specifically, the partition is
\[ A=\begin{bmatrix}  A_1 \\ \vdots \\ A_n \end{bmatrix}\in \mathbb{R}^{m\times m},\ B=\begin{bmatrix}  B_1,  \cdots , B_n \end{bmatrix}\in \mathbb{R}^{r\times r}, \]
with $ A_i\in \mathbb{R}^{m_i\times m}, B_i\in \mathbb{R}^{r\times r_i} $ and $ \sum_{i=1}^{n}m_i=m, \sum_{i=1}^{n}r_i=r. $
If there is a Lyapunov equation with $ B=A^T, $ it only needs to divide $ A $ by row, naturally. As for matrix $ C $, we can divide it like $ A $ by row or like $ B $ by column as required. The treatments for different partitions of $ C $ are similar, then we focus on the partition of $ C $ by column
\begin{equation}\label{C-division}
C=\begin{bmatrix}  C_1,  \cdots , C_n \end{bmatrix}\in \mathbb{R}^{m\times r}, \ C_i\in \mathbb{R}^{m\times r_i}.
\end{equation}

This paper aims to solve two main problems for the Sylvester equation:

\noindent(\uppercase\expandafter{\romannumeral1}) Least squares solution and exact solution case: Aim to get a least squares solution or an exact solution of the Sylvester equation \eqref{se} by solving a convex unconstrained optimization problem:
\begin{equation}\label{min:X_least}
\min_X||AX+XB-C||_F^2,
\end{equation}
which minimizes the squared Frobenius norm of residual matrix.

\noindent(\uppercase\expandafter{\romannumeral2}) Regularization case: Aim to solve the equation with penalty requirement to the unknown variable $ X $ by solving a regularized convex unconstrained optimization problem:
\begin{equation}\label{eq:penalty-centralized}
\min_{X} \frac{1}{2}||AX+XB-C||_F^2+\alpha g(X),
\end{equation}
where $ \alpha>0 $ is a tradeoff between minimizing the squared residual and the penalty term, $ g(\cdot) $ is a general convex function representing the regularization penalty of variables.

About the computation of matrix equations, an immediate idea is to use the Kronecker product to rewrite \eqref{se} as a standard linear algebraic equation
 $ (I_m\otimes A+B^T\otimes I_r)\mathrm{vec}(X)=\mathrm{vec}(C). $
It can be obtained that the solution of \eqref{se} for each $ \mathrm{vec}(C) $ is unique if and only if the matrix $ I_m\otimes A+B^T\otimes I_r $ is nonsingular, which is equivalent to requiring $ \mathrm{spec}(A)\cap \mathrm{spec}(-B)=\emptyset $ \cite{horn1991topics}. Different from centralized situations, the distributed situation may spoil the information structures because the sub-blocks of $ A $ and $ B $ are mixed up for the Kronecker product. Therefore, inspired by the work about distributed algorithms of the liner algebraic equation  \cite{shi2017networkTAC,liu2019arrow,liu2017asynchronous,wang2017double,cao2017continuous}, we have to study a new method, distinguished from the centralized and parallel algorithms, to solve the Sylvester equation \eqref{se} for given partition in a distributed way.

Consider a multi-agent network consisting of $ n $ agents, which is described by an undirected graph  $ \mathcal{G}(\mathcal{V}, \mathcal{E}, A_{\mathcal{G}}) $, where $ \mathcal{V}=\{1,\ldots, n\} $ is the set of nodes, $ \mathcal{E}\subset\mathcal{V}\times\mathcal{V} $ is the set of edges, and $ A_{\mathcal{G}}=[a_{i,j}]\in \mathbb{R}^{n\times n} $ is the adjacency matrix with $ a_{i,j}=a_{j,i}>0 $ if $ (i,j)\in \mathcal{E} $ and $ a_{i,j}=0 $ otherwise.
Denote the neighbor set of agent $ i $ as $ \mathcal{N}_i=\{j:(i,j)\in\mathcal{E}, j\in \mathcal{V}\}. $
The Laplacian matrix $ L_n=D-A_{\mathcal{G}} $, where $ D\in \mathbb{R}^{n\times n} $ is diagonal with $ D_{i,i}=\sum_{i=1}^{n}a_{i,j} $.
\begin{assumption}\label{ass:G-undirected-connected}
	The graph $ \mathcal{G} $ is undirected and connected.
\end{assumption}

Suppose that each agent knows a row sub-block of $ A $, a column sub-block of $ B $, and a column sub-block of $ C $, that is, agent $ i $ holds $ A_i, B_i, C_i $.
As a result, agent $ i $ in the network needs to estimate one common variable $ X^* $, which is a solution for the Sylvester equation and their final estimates are supposed to achieve consensus, that is, $ X_i=X^* $ for all $ i\in\mathcal{V} $.

Different from the cases given in \cite{zeng2018distributed}, to make $ AX $ and $ XB $ simultaneously distributed in design for the LRRC case, here we introduce two new variables $ Y, Z $ for coordinating the inconsistency due to the different block partitions of $ A $ and $ B $.  The following statement shows the idea.

For the LRRC case of matrices $ A $ and $ B $, and $ C $ divided as \eqref{C-division}, the Sylvester equation \eqref{se} is equivalent to
\begin{equation}\label{p1}
\left\{\begin{aligned}
&AX=Y\\&XB=C-Z\\&Y=Z
\end{aligned}\right.
\end{equation}
where supplementary variables $ Y $ and $ Z $ have corresponding block structures,
$ Y=\mathrm{col}\{Y_1, \ldots, Y_n\},
Z=[Z_1,  \cdots , Z_n] $
with $ Y_i\in \mathbb{R}^{m_i\times r}, Z_i\in \mathbb{R}^{m\times r_i} $.

In the considered multi-agent network, agent $i$ in the network $ \mathcal{G} $ has knowledge of information $ A_i, B_i, C_i $ and state set $ (X_i,Y_i, Z_i) $, and cooperates with its neighbors to compute $ X_i $ through the exchange of local information.

\begin{remark}
	If $ C $ is divided by row like $ A $, then the equivalent form of the Sylvester equation \eqref{se} is $ AX=C-Z,\ XB=Y,\ Z=Y $ where $ Z=\mathrm{col}\{Z_1, \ldots, Z_n\},$  $
	Y=\begin{bmatrix}  Y_1,  \cdots , Y_n \end{bmatrix} $ with $ Z_i\in \mathbb{R}^{m_i\times r}, Y_i\in \mathbb{R}^{m\times r_i} $. Algorithms and analysis ideas are analogous to the case in \eqref{p1}, so we will not discuss the case for space limitations.
\end{remark}

In order to reformulate the equivalent form \eqref{p1} as an optimization problem, the idea is to take one or more residuals of the equations as an objective function and take others as equality constrains. Unlike the unconstrained optimization problem \eqref{min:X_least}, the
reformulated problems are standard optimization problems with equality constraints.
Further, we need to consider the given distributed situation.
Since the coupling equality constraint $ Y=Z $ is not separable for each agent, and consensus constraints $ X_i=X_j $ for all $ i,j\in\mathcal{V} $ are hard to deal with directly, following transformations are necessary to assign each constraint to each agent clearly.

\begin{lemma}\label{lemma:constraints-eq}
	These two equivalent transformations hold for all $ i\in \mathcal{V} $.
	\begin{subequations}
		\begin{align}
		&Y=Z\ \Longleftrightarrow \exists \{W_i\}_{i=1}^n,~s.t.~ \nonumber\\ &\qquad ~~ [Y_i]_\rr-[Z_i]_\cc-\sum_{j=1}^{n}a_{i,j}(W_i-W_j)=0_{m\times r}, \label{eq:Y=Z}\\
		&X_i=X_j,\ \forall j\in \mathcal{V} \Longleftrightarrow \sum_{j=1}^{n}a_{ij}(X_i-X_j)=0_{m\times r}, \label{eq:X-i=X-j}
		\end{align}
	\end{subequations}
	where $ [a_{i,j}] $ is the corresponding adjacency matrix of an undirected and connected graph $\mathcal G$, $ W_i\in \mathbb{R}^{m\times r} $ are introduced to make up for the inconsistencies between $ [Y_i]_\rr $ and $ [Z_i]_\cc. $
\end{lemma}	

$ Proof $:
	As for \eqref{eq:Y=Z}, it is trivial that $ Y=Z $ holds by summing the right side of \eqref{eq:Y=Z} from $ i=1 $ to $ i=n $. On the other hand, according to the structures of $ Y $ and $ Z $, we have $ Y=\sum_{i=1}^{n}[Y_i]_\rr, Z=\sum_{i=1}^{n}[Z_i]_\cc. $
	Define $ L=L_n\otimes I_m, $ where $ L_n $ is the Laplacian matrix of graph $ \mathcal{G}. $ Therefore, $ \mathrm{ker}(L)=\{k\mathbf{1}_n\otimes I_m: k\in\mathbb{R}\} $ and
	$ \mathrm{Im}(L)=\ker(L)^{\bot}=\big\{w\in\mathbb{R}^{nm}:\sum_{i=1}^{n}w_i=0_m, w=\mathrm{col}\{w_1,\ldots,w_n\}, w_i\in\mathbb{R}^{m}, i\in \mathcal{V} \big\}. $
	Because $ \sum_{i=1}^{n}[Y_i]_\rr=\sum_{i=1}^{n}[Z_i]_\mathcal{C} $, the equation $ \sum_{i=1}^{n}([Y_i]_\rr-[Z_i]_\cc)=0_{m\times r} $  infers that every column of $ \mathrm{col}\{[Y_1]_\rr-[Z_1]_\cc, \cdots, [Y_n]_\rr-[Z_n]_\cc\} $ belongs to $ \mathrm{Im}(L). $ Thus, there exists $ W=\mathrm{col}\{W_1, \cdots, W_n\}\in \R^{nm\times r}, $ such that
	$LW=\mathrm{col}\{[Y_1]_\rr-[Z_1]_\cc, \cdots, [Y_n]_\rr-[Z_n]_\cc\}. $
	Then we have \eqref{eq:Y=Z}.
	
	As for \eqref{eq:X-i=X-j}, the equivalence is straightforward, and its proof is omitted here. \hfill$ \Box $
\begin{remark}
	The equality constraint \eqref{eq:Y=Z} is a
	typical coupling constraint, and \eqref{eq:X-i=X-j} is a
	typical  consensus constraint. The transformations in Lemma \ref{lemma:constraints-eq}  are always effective to deal with these two kinds of constraints.
\end{remark}

In brief, to reformulate the problem clearly, we have already introduced three types of variables, where their purposes are presented in Table \ref{table:variablesYZW}.
\begin{table*}[htpb]\small
	\centering
	\setlength{\abovecaptionskip}{0pt}
	\setlength{\belowcaptionskip}{0pt}
	\caption{
		The purposes of introduced variables}
	\label{table:variablesYZW}	
	\begin{tabular}{|c|c|}
		\hline
		\small Variable & Purpose \\
		\hline
		$ Y=\mathrm{col}\{Y_1,\cdots,Y_n\} $& $ Y_i=A_iX $, result of $ A_iX $ stored by sub-row of $ Y $\\
		\hline
		$ Z=[Z_1,\cdots,Z_n] $& $ Z_i=C_i-XB_i, $ result of $ C_i-XB_i $ stored by sub-column of $ Z $\\
		\hline
		$ W=\mathrm{col}\{W_1,\cdots,W_n\} $&  $ (L_n\otimes I_m)W=\mathrm{col}\{[Y_1]_\rr-[Z_1]_\cc, \cdots, [Y_n]_\rr-[Z_n]_\cc\} $ \\
		\hline
	\end{tabular}
\end{table*}
\section{The Cases of Least Squares Solution and Exact Solution}
\label{sec:exact+LS}
\subsection{Least Squares Solution Case}
\label{sec:least squares solution}
We first consider solving a least squares solution of the Sylvester equation \eqref{se} because it may not have exact solutions. Additionally, since local information is stored by different agents, it is not clear for each agent whether there are exact solutions. Of course, the proposed method for a least squares solution still works when there exist exact solutions.
\subsubsection{Optimization Model}
Reformulate \eqref{p1} as an optimization problem
\begin{equation}\label{eq:0-lsc}
\begin{aligned}
\min_{X,Y,Z}~ & \frac{1}{2}||XB-C+Z||_F^2\\
s.t.~
&AX=Y,\
Y=Z.
\end{aligned}
\end{equation}

Next, we distribute the information to every agent $ i $ and add the consensus constraints $ X_i=X_j $ for all $ i,j\in \mathcal{V} $. According to Lemma \ref{lemma:constraints-eq}, \eqref{eq:0-lsc} can be rewritten as a distributed optimization problem
\begin{equation}\label{lsc}
\begin{aligned}
\min_{\bar{X},Y,Z,W}~& \frac{1}{2}\sum_{i=1}^{n}||X_{i}B_i-C_i+Z_i||_F^2\\
s.t.~& \sum_{j=1}^{n}a_{i,j}(X_{i}-X_{j})=0_{m\times r}, \quad A_iX_{i}=Y_i,\\
&[Y_i]_\rr-[Z_i]_\cc-\sum_{j=1}^{n}a_{i,j}(W_i-W_j)=0_{m\times r}, \forall i\in \mathcal{V},
\end{aligned}
\end{equation}
with $ \bar{X}=\mathrm{col}\{X_1,\cdots, X_n\}, W=\mathrm{col}\{W_1,\cdots, W_n\} $ and $ Y, Z $ as in \eqref{p1}.
In the optimization problem \eqref{lsc}, $ \Lambda_i\in \mathbb{R}^{m\times r}, \Upsilon_i\in \mathbb{R}^{m_i\times r} $ and $ \Theta_i\in \mathbb{R}^{m\times r}, i\in \mathcal{V} $ are introduced as Lagrange multipliers associated with  these equality constraints and $ \Lambda=\mathrm{col}\{\Lambda_1,\cdots, \Lambda_n\}, \Upsilon=\mathrm{col}\{\Upsilon_1,\cdots, \Upsilon_n\}, \Theta=\mathrm{col}\{\Theta_1,\cdots, \Theta_n\} $.  Let $ (\bar{X}^*, Y^*, Z^*, W^*) $ and $ (\Lambda^*, \Upsilon^* , \Theta^*) $ be a pair of primal and dual optimal points. Then they satisfy the Karush-Kuhn-Tucker (KKT) condition, which leads to a necessary and sufficient condition for optimality, for all $ i\in\mathcal{V} $,
\begin{equation}\label{KKT-ls}
\left\{
\begin{array}{l}
0_{m\times r}=\sum_{j=1}^{n}a_{i,j}(X_{i}^*-X_{j}^*),\\
0_{m_i\times r}=A_iX_{i}^*-Y_i^*, \\
0_{m\times r}={[Y_i^*]}_\rr-[Z_i^*]_\cc-\sum_{j=1}^{n}a_{i,j}(W_i^*-W_j^*), \\
0_{m\times r}=-(X_{i}^*B_i-C_i+Z_i^*)B_i^T-\sum_{j=1}^{n}a_{i,j}(\Lambda_i^*-\Lambda_j^*)\\ \qquad~ \quad-A_i^T\Upsilon_i^*, \\
0_{m_i\times r}=\Upsilon_i^*-[I_{m_i}]_\cc\Theta_i^*,\\
0_{m\times r_i}=-(X_{i}^*B_i-C_i+Z_i^*)+\Theta_i^*[I_{r_i}]_\rr,\\
0_{m\times r}=\sum_{j=1}^na_{i,j}(\Theta_i^*-\Theta_j^*).
\end{array}\right.
\end{equation}

\subsubsection{Distributed Algorithm}
Each agent $ i $ in the network knows $ A_i, B_i, C_i $ and has its own state $ (X_i, Y_i, Z_i, W_i, \Lambda_i, \Upsilon_i, \Theta_i). $ Through the exchange of information, agent $ i $ can receive partial state information $ (X_j, W_j, \Theta_j), j\in \mathcal{N}_i $ from its neighbor agents. Then we propose Algorithm \ref{al:least-square} as the distributed algorithm for each agent $ i $ in the least squares solution case,
\begin{algorithm}\small
	\caption{Algorithm for the least squares solution case}
	\label{al:least-square}
	$ X_i(0)\in \R^{m\times r},  Y_i(0)\in \R^{m_i\times r}, Z_i(0)\in \R^{m\times r_i}, W_i(0)\in \R^{m\times r}, \Lambda_i(0)\in \R^{m\times r}, \Upsilon_i(0)\in \R^{m_i\times r}, \Theta_i(0)\in \R^{m\times r},  i \in \mathcal{V}, $
	\begin{equation}\label{sysal:least-square}
	\left\{
	\begin{array}{l}
	\dot{X}_i(t)=-(X_{i}(t)B_i-C_i+Z_i(t))B_i^T-A_i^T(A_iX_i(t)-Y_i(t))\\\qquad \qquad-A_i^T\Upsilon_i(t)-\sum_{j=1}^{n}a_{i,j}(\Lambda_i(t)-\Lambda_j(t))\\\qquad \qquad-\sum_{j=1}^{n}a_{i,j}(X_{i}(t)-X_{j}(t)),\\
	\dot{Y}_i(t)=\Upsilon_i(t)-[I_{m_i}]_\cc\Theta_i(t)+A_iX_i(t)-Y_i(t),\\
	\dot{Z}_i(t)=-(X_{i}(t)B_i-C_i+Z_i(t))+\Theta_i(t)[I_{r_i}]_\rr, \\
	\dot{W}_i(t)=\sum_{j=1}^{n}a_{i,j}(\Theta_i(t)-\Theta_j(t)), \\
	\dot{\Lambda}_i(t)=\sum_{j=1}^{n}a_{i,j}(X_{i}(t)-X_{j}(t)), \\
	\dot{\Upsilon}_i(t)=A_iX_i(t)-Y_i(t), \\
	\dot{\Theta}_i(t)=[Y_i(t)]_\rr-[Z_i(t)]_\mathcal{C}-\sum_{j=1}^{n}a_{i,j}(W_i(t)-W_j(t))\\\qquad \qquad-\sum_{j=1}^{n}a_{i,j}(\Theta_i(t)-\Theta_j(t)).
	\end{array}
	\right.
	\end{equation}
\end{algorithm}
where $ (X_i(t), Y_i(t), Z_i(t), W_i(t)) $ and $ (\Lambda_i(t), \Upsilon_i(t), \Theta_i(t)) $ are the estimates of solutions to problem \eqref{lsc} and the estimates of Lagrange multipliers by agent $ i $ at time $ t $, respectively. Algorithm \ref{al:least-square} can be viewed as the saddle-point dynamics of an augmented Lagrangian function $ \mathcal{L}_l $,
that is, for all $ i\in \mathcal{V}, $
\begin{align*}
&\dot{\varPhi}_i=-\nabla_{\varPhi_i}\mathcal{L}_l,~\text{when}~\varPhi_i=X_i,~Y_i,~Z_i~\text{or}~W_i,\\&
\dot{\varPsi}_i=\nabla_{\varPsi_i}\mathcal{L}_l,~\text{when}~\varPsi_i=\Lambda_i,~\Upsilon_i,~\text{or}~\Theta_i,
\end{align*}
where $ \mathcal{L}_l $ is defined in \eqref{eq:L_l}
and the last three terms are augmented terms.
\newcounter{mytempeqncnt}
	\begin{equation}\label{eq:L_l}
	\begin{aligned} \mathcal{L}_l=&\frac{1}{2}\sum_{i=1}^{n}||X_{i}B_i-C_i+Z_i||_F^2+\sum_{i=1}^{n}\langle\Lambda_i,\sum_{j=1}^{n}a_{i,j}(X_{i}-X_{j})\rangle_F+\sum_{i=1}^{n}\langle\Upsilon_i,A_iX_{i}-Y_i\rangle_F\notag\\
&+\sum_{i=1}^{n}\langle\Theta_i,[Y_i]_\rr-[Z_i]_\mathcal{C}-\sum_{j=1}^{n}a_{i,j}(W_i-W_j)\rangle_F-\frac{1}{2}\sum_{i=1}^{n}\langle\Theta_i,\sum_{j=1}^{n}a_{i,j}(\Theta_i-\Theta_j)\rangle_F\notag\\
&+\frac{1}{2}\sum_{i=1}^{n}\langle X_i,\sum_{j=1}^{n}a_{i,j}(X_i-X_j)\rangle_F+\frac{1}{2}\sum_{i=1}^{n}||A_iX_i-Y_i||_F^2
	\end{aligned}
	\end{equation}

\begin{lemma}\label{lemma:ls-solution-equal}
	Under Assumption \ref{ass:G-undirected-connected}, let $ P_l^*=\mathrm{col}\{\bar{X}^*,Y^*,Z^*,W^*,\Lambda^*,$  $\Upsilon^*,\Theta^*\}$. Then $
	(\bar{X}^*, Y^*, Z^*, W^*) $ is an optimal solution of problem \eqref{lsc} if and only if there exist $ \Lambda^*\in \mathbb{R}^{mn\times r}, \Upsilon^*\in \mathbb{R}^{m\times r}, $ and $ \Theta^*\in \mathbb{R}^{mn\times r} $, such that $ P_l^* $ is an equilibrium of system \eqref{sysal:least-square}.
\end{lemma}

This lemma holds obviously because of the KKT condition \eqref{KKT-ls}.

Denote $ P_l=\mathrm{col}\{\bar{X},Y,Z,W,\Lambda,\Upsilon,\Theta\}. $ The following result shows main convergence properties of Algorithm \ref{al:least-square}.
\begin{theorem}\label{thm:ls-convergence}
	Under Assumption \ref{ass:G-undirected-connected}, the following hold with Algorithm \ref{al:least-square}.
	\begin{enumerate}
		\item [(\romannumeral1)] Every equilibrium of the system \eqref{sysal:least-square} is Lyapunov stable and every solution is bounded;
		\item [(\romannumeral2)] The system \eqref{sysal:least-square} is globally semi-stable.
		Moreover, the solution $ P_l(t) $ with initial condition $ P_l(0) $ converges to an equilibrium of \eqref{sysal:least-square} exponentially;
		\item [(\romannumeral3)] For every agent $ i, $ the limit of its estimation $ \lim_{t\to \infty} X_i(t)=X^* $ is a least squares solution of the Sylvester equation \eqref{se}.
	\end{enumerate}
\end{theorem}

Its proof is in Subsection \ref{sub:proofs-ls-ex}.

\subsection{Exact Solution Case}
\label{sec:exact solution}
In the case when there exist exact solutions of the Sylvester equation \eqref{se}, we can provide a simplified algorithm with some conversions.
\begin{assumption}\label{ass:solution-exists}
	The Sylvester equation \eqref{se} has at least one exact solution.
\end{assumption}
\subsubsection{Optimization Model}
Let Assumptions \ref{ass:G-undirected-connected} and \ref{ass:solution-exists} hold. The set $ \{(X,Y,Z):AX=Y,XB=C-Z,Y=Z\} $ has at least one element. In this case, we can formulate the optimization problem \eqref{eq:0-lsc} as follows:

\begin{equation}\label{eq:0-esc}
\begin{aligned}
\min_{X,Y,Z}~& \frac{1}{2}||XB-C+Z||_F^2+\frac{1}{2}||AX-Y||_F^2\\
s.t.~&Y=Z.
\end{aligned}
\end{equation}
Equivalently, the problems \eqref{eq:0-lsc} and \eqref{eq:0-esc} have the same solution, because \eqref{eq:0-lsc}
can reach the optimal value $ 0. $
Rewrite the objective function as the sum of each agent's objective function; then $ \frac{1}{2}\sum_{i=1}^{n}||X_{i}B_i-C_i+Z_i||_F^2+\frac{1}{2}\sum_{i=1}^{n}||A_iX_{i}-Y_i||_F^2 $ can reach the minimum 0.
In fact,
\begin{equation} \label{dis-obj}
\begin{aligned}
&\frac{1}{2}\sum_{i=1}^{n}||X_{i}B_i-C_i+Z_i||_F^2+\frac{1}{2}\sum_{i=1}^{n}||A_iX_{i}-Y_i||_F^2\\
\le&  \frac{1}{2}\sum_{i=1}^{n}||X_{i}B_i-C_i+Z_i||_F^2+\frac{1}{2}\sum_{i=1}^{n}||A_iX_{i}-Y_i||_F^2\\&+\frac{1}{2}\sum_{i=1}^{n}\langle\sum_{j=1}^{n}a_{i,j}(X_{i}-X_{j}), X_i\rangle_F.
\end{aligned}\end{equation}
When the right-hand side term reaches 0, $ X_i=X_j $ must hold for all $ i,j\in \mathcal{V} $. If the left-hand side term reaches 0 with $ X_i=X_j $ as an exact solution of the original equation, the right-hand side term equals 0 as well. Clearly, we can use the right-hand side of  \eqref{dis-obj} as the objective function to replace consensus constraints $ X_i=X_j. $ Then the following distributed optimization problem is formulated:
\begin{equation}\label{esc}
\begin{aligned}
\min_{\bar{X},Y,Z,W}&~ \frac{1}{2}\sum_{i=1}^{n}||X_{i}B_i-C_i+Z_i||_F^2+\frac{1}{2}\sum_{i=1}^{n}||A_iX_{i}-Y_i||_F^2\\&+\frac{1}{2}\sum_{i=1}^{n}\langle\sum_{j=1}^{n}a_{i,j}(X_{i}-X_{j}), X_i\rangle_F\\
s.t.
&~[Y_i]_\rr-[Z_i]_\cc-\sum_{j=1}^{n}a_{i,j}(W_i-W_j)=0_{m\times r}, \forall i\in \mathcal{V}.
\end{aligned}
\end{equation}

\begin{remark}
	Denote  an optimal point for problem \eqref{esc} by $ (\bar{X}^*, Y^*, Z^*,W^*) $.
	According to Assumption \ref{ass:solution-exists}, the objective function can reach 0, that is, $ X_i^*B_i-C_i+Z_i^*=0, A_iX_i^*-Y_i^*=0 $ and $ X_i^*=X_j^*, $ where $ X_i^*=X_j^*=X^* $ is an exact solution of \eqref{se} and $ (X^*, Y^*, Z^*) $ is an exact solution of \eqref{p1}.
\end{remark}

In the optimization problem \eqref{esc}, $ \Theta_i\in \mathbb{R}^{m\times r} $ is introduced as the Lagrange multiplier associated with the $ i$-th equality constraint and $ \Theta=\mathrm{col}\{\Theta_1,\cdots, \Theta_n\}. $  Let $ (\bar{X}^*, Y^*, Z^*, W^*) $ and $ \Theta^* $ be a pair of primal and dual optimal points. Then they satisfy the corresponding KKT condition \eqref{KKT-es}, which proposes necessary and sufficient conditions for the optimality, for all $ i\in\mathcal{V}, $
\begin{equation}\label{KKT-es}
\left\{
\begin{array}{l}
0_{m\times r}={[Y_i^*]}_\rr-[Z_i^*]_\cc-\sum_{j=1}^{n}a_{i,j}(W_i^*-W_j^*), \\
0_{m\times r}=-(X_{i}^*B_i-C_i+Z_i^*)B_i^T-\sum_{j=1}^{n}a_{i,j}(\Lambda_i^*-\Lambda_j^*)\\\qquad~~~~-\sum_{j=1}^na_{i,j}(X^*_i-X^*_j), \\
0_{m_i\times r}=(A_iX_i^*-Y^*_i)-[I_{m_i}]_\cc\Theta_i^*,\\
0_{m\times r_i}=-(X_{i}^*B_i-C_i+Z_i^*)+\Theta_i^*[I_{r_i}]_\rr,\\
0_{m\times r}=\sum_{j=1}^na_{i,j}(\Theta_i^*-\Theta_j^*).
\end{array}\right.
\end{equation}

\subsubsection{Distributed Algorithm}
Each agent $ i $ in the network knows $ A_i, B_i, C_i $ and has its own state $ (X_i, Y_i, Z_i, W_i, \Theta_i). $ Through the exchange of information, agent $ i $ can receive partial state information $ (X_j, W_j, \Theta_j), j\in \mathcal{N}_i $ from its neighbor agents. Then we propose Algorithm \ref{al:exact-solution}, using the similar idea of Algorithm \ref{al:least-square}, as the distributed algorithm for each agent $ i $ in the exact solution case,
\begin{algorithm}[h]\small
	\caption{Algorithm for the exact solution case}
	\label{al:exact-solution}
	$ X_i(0)\in \R^{m\times r},  Y_i(0)\in \R^{m_i\times r}, Z_i(0)\in \R^{m\times r_i}, W_i(0)\in \R^{m\times r}, \Theta_i(0)\in \R^{m\times r},  i\in \mathcal{V}, $
	\begin{equation}\label{sysal:exact-solution}\left\{
	\begin{array}{l}
	\dot{X}_i(t)=-\left(X_{i}(t)B_i-C_i+Z_i(t)\right)B_i^T-A_i^T\left(A_iX_i(t)-Y_i(t)\right)\\\qquad \qquad-\sum_{j=1}^{n}a_{i,j}(X_{i}(t)-X_{j}(t)),\\
	\dot{Y}_i(t)=\left(A_iX_i(t)-Y_i(t)\right)-[I_{m_i}]_\cc\Theta_i(t), \\
	\dot{Z}_i(t)=-\left(X_{i}(t)B_i-C_i+Z_i(t)\right)+\Theta_i(t)[I_{r_i}]_\rr, \\
	\dot{W}_i(t)=\sum_{j=1}^{n}a_{i,j}(\Theta_i(t)-\Theta_j(t)), \\
	\dot{\Theta}_i(t)=[Y_i(t)]_\rr-[Z_i(t)]_\cc-\sum_{j=1}^{n}a_{i,j}(W_i(t)-W_j(t))\\\qquad \qquad-\sum_{j=1}^{n}a_{i,j}(\Theta_i(t)-\Theta_j(t)).
	\end{array}\right.
	\end{equation}
\end{algorithm}
where $ (X_i(t), Y_i(t), Z_i(t), W_i(t)) $ and $ \Theta_i(t) $ are the estimates of solutions to problem \eqref{esc} and the estimate of Lagrange multiplier by agent $ i $ at time $ t $, respectively.
\begin{remark}
	Compared with Algorithm \ref{al:least-square}, the state variables that every agent needs to calculate in Algorithm \ref{al:exact-solution} have changed from $ (X_i(t), Y_i(t), Z_i(t), W_i(t), \Lambda_i(t), $  $\Upsilon_i(t), \Theta_i(t)) $ to $ (X_i(t), Y_i(t), Z_i(t), \\ W_i(t), \Theta_i(t)) $.
	On the premise of exact solutions, the algorithm dimension has been reduced to a certain extent.	
\end{remark}

\begin{lemma}\label{lemma:exactcase-solution-equal}
	Under Assumptions \ref{ass:G-undirected-connected} and \ref{ass:solution-exists}, let $ P_e^*=\mathrm{col}\{\bar{X}^*, Y^*,Z^*,$  $W^*,\Theta^*\} $. Then $ (\bar{X}^*, Y^*, Z^*, W^*) $ is an optimal solution of problem \eqref{esc} if and only if there exists $ \Theta^*\in \mathbb{R}^{mn\times r} $ such that $ P_e^* $ is an equilibrium of system \eqref{sysal:exact-solution}.
\end{lemma}

This lemma can be easily proved by the KKT optimality condition \eqref{KKT-es}.

Denote $ P_e=\mathrm{col}\{\bar{X},Y,Z,W,\Theta\} $. The following result shows main convergence properties of Algorithm \ref{al:exact-solution}.

\begin{theorem}\label{thm:exact-converge}
	Under Assumptions \ref{ass:G-undirected-connected} and  \ref{ass:solution-exists}, the following hold with Algorithm \ref{al:exact-solution}.
	\begin{enumerate}
		\item[(\romannumeral1)] Every equilibrium of the system \eqref{sysal:exact-solution} is Lyapunov stable and every solution is bounded;
		\item[(\romannumeral2)] The system \eqref{sysal:exact-solution} is globally semi-stable. Moreover, the solution $ P_e(t) $ with initial condition $ P_e(0) $ converges to an equilibrium of \eqref{sysal:exact-solution} exponentially;
		\item[(\romannumeral3)] For every agent $ i $, the limit of its estimation $ \lim_{t\to \infty} X_i(t)=X^* $ is an exact solution of the Sylvester equation \eqref{se}.
	\end{enumerate}
\end{theorem}

Its explanation is in the next subsection.

\subsection{Proofs of Theorems \ref{thm:ls-convergence} and \ref{thm:exact-converge}}\label{sub:proofs-ls-ex}
	$ The~Proof~of~Theorem~\ref{thm:ls-convergence} $:
	(\romannumeral1) First, note that $ P_l^*=\mathrm{col} \{\bar{X}^*, Y^*, Z^*, W^*, \Lambda^*, \Upsilon^*,$  $ \Theta^*\} $ is an equilibrium of \eqref{sysal:least-square}, which satisfies the KKT condition \eqref{KKT-ls}. Define a positive definite function:
	$ V_l(P_l):=\frac{1}{2}||P_l-P_l^*||_F^2. $
	By using Algorithm \ref{al:least-square} and the KKT condition \eqref{KKT-ls},
	we calculate the derivative of $ V_l $ with respect to time $ t $ and obtain
	\begin{align*}
	\frac{d}{dt}V_l=&-\sum_{i=1}^{n}||(X_i-X_i^*)B_i+(Z_i-Z_i^*)||_F^2 \\&-\sum_{i=1}^{n}||A_i(X_i-X_i^*)-(Y_i-Y_i^*)||_F^2\\&-\frac{1}{2}\sum_{i=1}^{n}\sum_{j=1}^{n}a_{i,j}(||X_i-X_j||_F^2+||\Theta_i-\Theta_j||_F^2)\le 0.
	\end{align*}
	Because $ V_l $ is positive definite and $ \frac{d}{dt}V_l\le 0, $ the equilibrium $ P_l^* $ of \eqref{sysal:least-square} is Lyapunov stable. Moreover, because of $ V_l\to \infty $ as $ ||P_l-P_l^*||_F\to \infty $, each solution with given initial conditions of \eqref{sysal:least-square} is bounded for all $ t\ge 0 $ (Theorem 4 in \cite{laSalle1960some}).
	
	\noindent(\romannumeral2) With the condition of $ \frac{d}{dt}V_l=0, $ we define the set
	$ \rr_l=\{P_l:\frac{d}{dt}V_l=0\}
	=\{P_l: X_i=X_j, \Theta_i=\Theta_j, (X_i-X_i^*)B_i=-(Z_i-Z_i^*), A_i(X_i-X_i^*)=Y_i-Y_i^*\}. $	
	Let $ \mathcal{M}_l $ be the largest invariant subset of $ \mathcal{\bar{R}}_l. $ It follows from the LaSalle's invariance principle \cite{laSalle1960some} that
	$ P_l(t)=\mathrm{col}\{\bar{X}(t), Y(t), Z(t), W(t), \Lambda(t), \Upsilon(t), \Theta(t)\}\to \mathcal{M}_l,$ as $ t\to \infty $
	and $ \mathcal{M}_l $ is positively invariant.
	
	Assume that $ \hat{P_l}(t)\in \mathcal{M}_l $ for all $ t\ge 0 $ is a trajectory of \eqref{sysal:least-square}.
	For all $ i,j\in \mathcal{V} $, $ \hat{P_l}(0)\in \mathcal{M}_l $,
	$ \hat{X}_i(t)=\hat{X}_j(t), \hat{\Theta}_i(t)=\hat{\Theta}_j(t), (\hat{X}_i(t)-\hat{X}_i^*)B_i=-(\hat{Z}_i(t)-\hat{Z}_i^*), $ and $ A_i(\hat{X}_i(t)-\hat{X}_i^*)=\hat{Y}_i(t)-\hat{Y}_i^*. $ Then we have
	\begin{itemize}
		\item [a)]$ \dot{\hat{\Lambda}}_i(t)=0_{m\times r} $ because $ \hat{X}_i(t)=\hat{X}_j(t).$
		\item [b)]$ \dot{\hat{W}}_i(t)=0_{m\times r} $ because  $ \hat{\Theta}_i(t)=\hat{\Theta}_j(t).$
		\item [c)] Since $ A_i(\hat{X}_i(t)-\hat{X}_i^*)=\hat{Y}_i(t)-\hat{Y}_i^* $ and $ A_i\hat{X}_i^*=\hat{Y}_i^* $, there hold $ A_i\hat{X}_i(t)=\hat{Y}_i(t)$ and $ \dot{\hat{\Upsilon}}_i(t)=0_{m_i\times r}.$
		\item   [d)]$\dot{\hat{X}}_i(t)=-(\hat{X}_{i}(t)B_i-C_i+\hat{Z}_i(t))B_i^T-\sum_{j=1}^{n}a_{i,j}(\hat{\Lambda}_i(t)-\hat{\Lambda}_j(t))-A_i^T\hat{\Upsilon}_i(t)-\sum_{j=1}^{n}a_{i,j}(\hat{X}_{i}(t)-\hat{X}_{j}(t))-A_i^T(A_i\hat{X}_i(t)-\hat{Y}_i(t))
		=-(\hat{X}_i^*B_i-C_i+\hat{Z}_i^*)B_i^T-\sum_{j=1}^{n}(\hat{\Lambda}_i(0)-\hat{\Lambda}_j(0))-A_i^T\hat{\Upsilon}_i(0) $	
		is a constant matrix. If $ \dot{\hat{X}}_i(t)\neq0_{m\times r}, \hat{X}_i(t) $ will be unbound, which is a contradiction from the boundedness in (\romannumeral1), so $ \dot{\hat{X}}_i(t)=0_{m\times r} $.
		\item [e)]  Note that
		$ \hat{Y}_i(t)=A_i\hat{X}_i(t)=A_i\hat{X}_i(0) $ is a constant matrix. Then $ \dot{\hat{Y}}_i(t)=0_{m_i\times r}$.
		\item [f)]  Due to $ (\hat{X}_i(t)-\hat{X}_i^*)B_i=-(\hat{Z}_i(t)-\hat{Z}_i^*), $ and $ \dot{\hat{X}}_i(t)=0_{m\times r},  \hat{Z}_i(t)=\hat{Z}_i^*-(\hat{X}_i(0)-\hat{X}_i^*)B_i $ is a constant matrix. Because of the boundedness of the solution in (\romannumeral1), $ \dot{\hat{Z}}_i(t)=0_{m\times r_i}. $		
		\item [g)]  Recall that $ \dot{\hat{\Theta}}_i(t) $ in \eqref{sysal:least-square} is a constant matrix because $ \dot{\hat{Y}}_i(t)=0_{m_i\times r}, \dot{\hat{Z}}_i(t)=0_{m\times r_i}, \dot{\hat{W}}_i(t)=0_{m\times r_i} $. Then $ \dot{\hat{\Theta}}_i(t)=0_{m\times r}.$
	\end{itemize}
	To sum up, 	$ \mathcal{M}_l\subset\{P_l(t): \dot{X}_i(t)=0_{m\times r}, \dot{Y}_i(t)=0_{m_i\times r},$  $ \dot{Z}_i(t)=0_{m\times r_i}, \dot{W}_i(t)=0_{m\times r}, \dot{\Lambda}_i(t)=0_{m\times r}, \dot{\Upsilon}_i(t)=0_{m_i\times r}, \dot{\Theta}_i(t)=0_{m\times r} \} $ and any point in $ \mathcal{M}_l $ is an equilibrium of \eqref{sysal:least-square}, which is Lyapunov stable. In view of Lemma \ref{lemma:converge} and the fact that every solution is bounded, \eqref{sysal:least-square} is globally semi-stable. In other words, for any initial conditions $ P_l(0) $, with Algorithm \ref{al:least-square} the solution $ P_l(t) $ converges to a Lyapunov stable equilibrium of \eqref{sysal:least-square}.
	
	If, for any given initial conditions, the solution of a linear time-invariant system converges, then it converges exponentially \cite{zeng2018distributed}.
	In fact, \eqref{sysal:least-square} is a linear time-invariant system, hence, every solution for
	given initial conditions of system \eqref{sysal:least-square} converges to an equilibrium exponentially.
	
	\noindent(\romannumeral3) Evidently, due to Lemma \ref{lemma:ls-solution-equal}, for given initial value $ P_l(0) $, every trajectory $ (\bar{X}(t), Y(t), Z(t), W(t)) $  converges to an equilibrium
	which is an optimal solution of problem \eqref{lsc}. Thus, $ X_i^*=X_j^*, $ for all $ i,j\in \mathcal{V} $ from the equality constraints. In
	other word, for every agent $ i $, the limit of its estimation $ \lim_{t\to \infty}X_i(t)=X^* $ is a least squares solution of the Sylvester equation \eqref{se}. \hfill $ \Box $
	
	As for Theorem \ref{thm:exact-converge}, we define a similar Lyapunov function   $ V_e(P_e):=\frac{1}{2}||P_e-P_e^*||_F^2, $  whose derivative with respect to time $ t $ satisfying $ \frac{d}{dt}V_e\le 0. $
	Then Theorem \ref{thm:exact-converge} can be proved by using the analogical analytical method about the proof of Theorem \ref{thm:ls-convergence}. More details are omitted here.	
	
	\section{Regularization Case}
	\label{sec:regularization}
	The regularization technique, such as Tikhonov regularization and LASSO \cite{tibshirani1996regression}, has been widely applied to statistics and machine learning to deal with  practical problems.
	Similar ideas can also be applied to the study of matrices, whose regularization also has many applications, such as multi-task learning, matrix completion and multivariate regression. Matrix regularization penalty functions are typically chosen to be convex.  Specifically, the $ l_1 $-norm can be selected to enforce sparsity; the $ l_{2,1} $-norm can be selected to enforce structured sparsity in feature selection; and the nuclear norm can be selected to enforce low rank approximation properties \cite{yin2008bregman,nie2010efficient,Candes2009exact,recht2010guaranteed,li2017constrained}.
	
	For the problem \eqref{eq:penalty-centralized}, we present the following regularization distributed optimization problem with a constant $ \alpha>0, $
	\begin{equation}\label{nco}
	\begin{aligned}
	\min_{\bar{X},Y,Z,W}~ & \frac{1}{2}\sum_{i=1}^{n}||X_{i}B_i-C_i+Z_i||_F^2+\alpha\sum_{i=1}^{n}g(X_i)\\
	s.t.~ & \sum_{j=1}^{n}a_{i,j}(X_{i}-X_{j})=0_{m\times r}, \quad A_iX_{i}=Y_i,\\
	&[Y_i]_\rr-[Z_i]_\cc-\sum_{j=1}^{n}a_{i,j}(W_i-W_j)=0_{m\times r},\forall i\in \mathcal{V}.
	\end{aligned}
	\end{equation}
	Then we introduce $ \Lambda_i\in \mathbb{R}^{m\times r}, \Upsilon_i\in \mathbb{R}^{m_i\times r} $ and $ \Theta_i\in \mathbb{R}^{m\times r}, i\in \mathcal{V} $ as the Lagrange multipliers associated with  these equality constraints and $ \Lambda=\mathrm{col}\{\Lambda_1,\cdots, \Lambda_n\},$  $ \Upsilon=\mathrm{col}\{\Upsilon_1,\cdots, \Upsilon_n\},$  $ \Theta=\mathrm{col}\{\Theta_1,\cdots, \Theta_n\} $.  Let $ (\bar{X}^*, Y^*, Z^*, W^*) $ and $ (\Lambda^*, \Upsilon^* , \Theta^*) $ be a pair of primal and dual optimal points. They satisfy the KKT condition, which leads to a necessary and sufficient condition for optimality, for all $ i\in\mathcal{V} $,
	\begin{equation}\label{KKT-nco}
	\left\{
	\begin{array}{l}
	0_{m\times r}=\sum_{j=1}^{n}a_{i,j}(X_{i}^*-X_{j}^*),\\
	0_{m_i\times r}=A_iX_{i}^*-Y_i^*, \\
	0_{m\times r}={[Y_i^*]}_\rr-[Z_i^*]_\cc-\sum_{j=1}^{n}a_{i,j}(W_i^*-W_j^*), \\
	0_{m\times r}\in -(X_i^*B_i-C_i+Z_i^*)B_i^T-\alpha\partial g(X_i)\\\qquad~~~~-\sum_{j=1}^{n}a_{i,j}(\Lambda_i^*-\Lambda_j^*)-A_i^T\Upsilon_i^*,\\
	0_{m_i\times r}=\Upsilon_i^*-[I_{m_i}]_\cc\Theta_i^*+A_iX_i^*-Y_i^*,\\
	0_{m\times r_i}=-(X_{i}^*B_i-C_i+Z_i^*)+\Theta_i^*[I_{r_i}]_\rr,\\
	0_{m\times r}=\sum_{j=1}^{n}a_{i,j}(\Theta_i^*-\Theta_j^*).
	\end{array}\right.
	\end{equation}
	The fourth formula in \eqref{KKT-nco} means that there exists $ h_i^* \in \partial g(X_i^*) $ such that
	\begin{equation}\label{eq:hi^*}
	\begin{aligned}
	0_{m\times r}=&-(X_i^*B_i-C_i+Z_i^*)B_i^T-\alpha h_i^*\\&-\sum_{j=1}^{n}a_{i,j}(\Lambda_i^*-\Lambda_j^*)-A_i^T\Upsilon_i^*.
	\end{aligned}\end{equation}
	Note that each agent $ i $ in the network can get access to the information about $ A_i, B_i, C_i $ and has its own state $ (\dot{X_i}, X_i, Y_i, Z_i, W_i, \Lambda_i, \Upsilon_i, \Theta_i). $ Through the exchange of information, agent $ i $ can receive partial state information $ (\dot{X_j}, X_j, W_j, \Theta_j), $  $ j\in \mathcal{N}_i $ from its neighbor agents. Then we propose Algorithm \ref{al:regular} as the distributed algorithm for each agent $ i $ in the regularization case,
	\begin{algorithm}\small
		\caption{Algorithm for the regularization case}
		\label{al:regular}
		$ X_i(0)\in \R^{m\times r},  Y_i(0)\in \R^{m_i\times r}, Z_i(0)\in \R^{m\times r_i}, W_i(0)\in \R^{m\times r}, \Lambda_i(0)\in \R^{m\times r}, \Upsilon_i(0)\in \R^{m_i\times r}, \Theta_i(0)\in \R^{m\times r},   i\in \mathcal{V}, $
		\begin{equation}\label{sysal:regular}
		\left\{
		\begin{array}{l}
		\dot{X}_i(t)\in-\left(X_{i}(t)B_i-C_i+Z_i(t)\right)B_i^T-\alpha\partial g(X_i)-A_i^T\Upsilon_i(t)\\\qquad \qquad-\sum_{j=1}^{n}a_{i,j}(\Lambda_i(t)-\Lambda_j(t))- A_i^T\left(A_iX_i(t)-Y_i(t)\right)\\\qquad \qquad-\sum_{j=1}^{n}a_{i,j}(X_{i}(t)-X_{j}(t)) ,\\
		\dot{Y}_i(t)=\Upsilon_i(t)-[I_{m_i}]_\cc\Theta_i(t)+A_i\left(X_i(t)-\dot{X_i}(t)\right)-Y_i(t),\\
		\dot{Z}_i(t)=-\left(\left(X_{i}(t)-\dot{X_i}(t)\right)B_i-C_i+Z_i(t)\right)+\Theta_i(t)[I_{r_i}]_\rr,\\
		\dot{W}_i(t)=\sum_{j=1}^{n}a_{i,j}\left(\Theta_i(t)-\Theta_j(t)\right),\\
		\dot{\Lambda}_i(t)=\sum_{j=1}^{n}a_{i,j}\left(X_{i}(t)-X_{j}(t)+\dot{X_i}(t)-\dot{X_j}(t)\right),\\
		\dot{\Upsilon}_i(t)=A_i\left(X_i(t)+\dot{X_i}(t)\right)-Y_i(t),\\
		\dot{\Theta}_i(t)=[Y_i(t)]_\rr-[Z_i(t)]_\cc-\sum_{j=1}^{n}a_{i,j}\left(W_i(t)-W_j(t)\right)\\\qquad \qquad-\sum_{j=1}^{n}a_{i,j}(\Theta_i(t)-\Theta_j(t)).
		\end{array}\right.
		\end{equation}
	\end{algorithm}
	where $ (X_i(t), Y_i(t), Z_i(t), W_i(t)) $ and $ (\Lambda_i(t), \Upsilon_i(t), \Theta_i(t)) $ are estimates of solutions to problem \eqref{nco} and estimates of Lagrange multipliers by agent $ i $ at time $ t $, respectively.  Algorithm
	\ref{al:regular} is inspired by the saddle-point dynamics and derivative feedback. Derivative feedback is a useful method to make the algorithm convergent \cite{antipin1994feedback,zeng2018disfornonsmooth}.
	
	\begin{lemma}\label{lemma:non-anops=an-equilibrium}
		Under  Assumption \ref{ass:G-undirected-connected}, let $ P_r^*=\mathrm{col}\{\bar{X}^*, Y^*, Z^*, W^*, \Lambda^*,$  $ \Upsilon^*, \Theta^*\}.$ Then $ (\bar{X}^*, Y^*, Z^*, W^*) $ is an optimal solution of problem \eqref{nco} if and only if there exist $ \Lambda^*\in \mathbb{R}^{mn\times r}, \Upsilon^*\in \mathbb{R}^{m\times r}, $ and $ \Theta^*\in \mathbb{R}^{mn\times r} $, such that $ P_r^* $ is an equilibrium of system \eqref{sysal:regular}.
	\end{lemma}
	
	Obviously, Lemma \ref{lemma:non-anops=an-equilibrium} can be easily proved by the KKT condition \eqref{KKT-nco}.
	
	Denote  $\dot{X_i}(t)\in \mathcal{F}_{X_i}, P_r=\mathrm{col}\{\bar{X},Y,Z,W,\Lambda,\Upsilon,\Theta\} $  and $ \dot{P_r}\in \mathcal{F}_{P_r}(P_r(t)). $
	Let $ P_r^* $ be an equilibrium of system \eqref{sysal:regular}.
	The next result shows main convergence properties of Algorithm \ref{al:regular}.
	
	\begin{theorem}\label{thm:regular-converge}
		Under Assumption \ref{ass:G-undirected-connected},
		the following hold with Algorithm \ref{al:regular}.
		\begin{enumerate}
			\item[(\romannumeral1)] Every equilibrium of the system \eqref{sysal:regular} is Lyapunov stable and every solution is bounded;
			\item[(\romannumeral2)] The system \eqref{sysal:regular} is globally semi-stable.
			\item[(\romannumeral3)] For every agent $ i, $ the limit of its estimation $ \lim_{t\to \infty} X_i(t)=X^* $ is an optimal solution of the problem \eqref{eq:penalty-centralized} based on the Sylvester equation \eqref{se}.
		\end{enumerate}
	\end{theorem}
	
	The proof of Theorem \ref{thm:regular-converge} follows some statements.
	According to Lemma \ref{lemma:df-dt} and the convexity of $ g $, $ \frac{d}{dt}g(X_i(t))\in \tilde{\mathcal{L}}_{\mathcal{F}_{X_i}}g, $ that is $ \frac{d}{dt}g(X_i(t))=\langle h_i,\dot{X_i}\rangle_F, $ for all $ h_i\in \partial g(X_i). $
	Because $ g(\cdot) $ is a convex function, for any $ h_i^*\in \partial g(X_i^*) $ and any $ h_i\in \partial g(X_i), $ there holds
	\begin{equation}\label{eq:h_i-h_i^*}
	\langle h_i-h_i^*,X_i-X_i^*\rangle_F\ge 0, \quad \forall  i\in \mathcal{V}.
	\end{equation}
	
	Define a scalar function \begin{align*}
	V_1:=&\sum_{i=1}^{r}\left(\alpha g(X_i)-\alpha g(X_i^*)-\alpha\langle h_i^*,X_i-X_i^*\rangle_F\right)\\&+\frac{1}{4}\sum_{i=1}^{r}\sum_{j=1}^{r}a_{i,j}||X_i-X_j||_F^2,
	\end{align*}
	where $ h_i^* $ satisfies \eqref{eq:hi^*}. Then $ V_1 $ is nonnegative.	
	Define a scalar function $ V_r(P_r):= V_1+V_2, $ where
	\begin{align*}
	V_2:=&\frac{1}{2}||P_{r}-P_{r}^*||_F^2+\frac{1}{2}\sum_{i=1}^{r}||A_i(X_i-X_i^*)||_F^2\\&+\frac{1}{2}\sum_{i=1}^{r}||(X_i-X_i^*)B_i||_F^2.
	\end{align*}
	Then the function $ V_r $ has the following properties:
	\begin{enumerate}
		\item [(a)] $ V_r(P_{r}^*)=0, V_r(P_{r})>0 $ for all $ P_{r}\ne P_{r}^*; $ $ V\to \infty $ as $ ||P_r-P_r^*||_F\to \infty $.
		\item [(b)] $ \frac{d}{dt}V_r(P_r(t))\in \tilde{\mathcal{L}}_{\mathcal{F}_{P_r}}V_r(P_r)=\{\nabla V^Ty:y\in\mathcal{F}_{P_r}(P_r)\}. $
	\end{enumerate}
	
	Next, it is time to prove Theorem \ref{thm:regular-converge}.

$ The~Proof~of~ Theorem~\ref{thm:regular-converge} $:
		(\romannumeral1) Note that $ P_r^*=\mathrm{col}\{\bar{X}^*,Y^*,Z^*,W^*,\Lambda^*,\Upsilon^*,\Theta^*\} $ is an equilibrium of \eqref{sysal:regular}. First, we calculate $ \frac{d}{dt}V_1(P_r(t)), $ for $ h_i\in \partial g(X_i), $ denote $ \dot{X_i}=-\left(X_{i}B_i-F_i+Z_i\right)B_i^T-\alpha H_i(X_i)-\sum_{j=1}^{n}a_{i,j}(\Lambda_i-\Lambda_j)-A_i^T\Upsilon_i-\sum_{j=1}^{n}a_{i,j}(X_{i}-X_{j})- A_i^T\left(A_iX_i-Y_i\right), $ then
		\begin{align*}
		&\frac{d}{dt} {V_1}=\sum_{i=1}^{r}\langle \alpha H_i(X_i)-\alpha h_i^*+\sum_{j=1}^{r}a_{i,j}(X_i-X_j),\dot{X_i}\rangle_F\\
		=&-\sum_{i=1}^{r}\langle \dot{X_i}, \dot{X_i}\rangle_F-\sum_{i=1}^{r}\langle \left(X_{i}B_i-F_i+Z_i\right)B_i^T+\sum_{j=1}^{r}a_{i,j}(\Lambda_i\\&-\Lambda_j)+A_i^T\Upsilon_i+A_i^T\left(A_iX_i-Y_i\right)+\alpha h_i^*,\dot{X_i}\rangle_F.
		\end{align*}
		Denote the second term in $ \frac{d}{dt} {V_1} $ by $ \beta$  
		With the help of the definition of $ h_i^* $, $ \beta $ can be rewritten as
		$ \beta=
		-\sum_{i=1}^{n}\langle (X_i-X_i^*)B_iB_i^T, \dot{X_i}\rangle_F-\sum_{i=1}^{n}\langle (Z_i-Z_i^*)B_i^T, \dot{X_i}\rangle_F-\sum_{i=1}^{n}\sum_{j=1}^{n}a_{i,j}\langle \Lambda_i-\Lambda_i^*, \dot{X_i}-\dot{X_j}\rangle_F-\sum_{i=1}^{n}\langle A_i^T(\Upsilon_i-\Upsilon_i^*), \dot{X_i}\rangle_F-\sum_{i=1}^{n}\langle A_i^T(A_i(X_i-X_i^*), \dot{X_i}\rangle_F+\sum_{i=1}^{n}\langle A_i^T(Y_i-Y_i^*), \dot{X_i}\rangle_F. $	
		Next, we calculate $ \frac{d}{dt}V_2. $
		Then, from $ V_r=V_1+V_2, $ we have the derivative of $ V_r $ with respect to time $ t $,
		\begin{align*}
		&\frac{d}{dt}V_r= -\sum_{i=1}^{n}||(X_i-X_i^*)B_i+(Z_i-Z_i^*)||_F^2\\&-\sum_{i=1}^{n}||A_i(X_i-X_i^*)-(Y_i-Y_i^*)||_F^2\\&-\frac{1}{2}\sum_{i=1}^{n}\sum_{j=1}^{n}a_{i,j}||X_i-X_j||_F^2-\frac{1}{2}\sum_{i=1}^{n}\sum_{j=1}^{n}a_{i,j}||\Theta_i-\Theta_j||_F^2\\&-\sum_{i=1}^{n}||\dot{X_i}||_F^2-\alpha\langle X_i-X_i^*, h_i-h_i^*\rangle_F\le 0.
		\end{align*}
		Thus, every element in $ \tilde{\mathcal{L}}_{\mathcal{F}_{P_r}}V_r $ is nonpositive; that is $ \max\tilde{\mathcal{L}}_{\mathcal{F}_{P_r}}V\le 0. $ As a result, from the properties of $ V_r $ and Lemma \ref{lemma:nonsmooth_stable}, the equilibrium $ P^*_r $ is Lyapunov stable.
		Moreover, because $ V_r $ is radially unbounded, each solution with given initial conditions of system \eqref{sysal:regular} is bounded for all $ t\ge 0 $.
		
		\noindent(\romannumeral2)
		Define the set
		$ \mathcal{R}_r=\{P_r:\frac{d}{dt}V_r=0\}
		=\{P_r: \dot{X_i}=0, X_i=X_j, \Theta_i=\Theta_j,  (X_i-X_i^*)B_i=-(Z_i-Z_i^*),\ \ A_i(X_i-X_i^*)=Y_i-Y_i^*,\langle X_i-X_i^*, h_i-h_i^*\rangle_F=0\}\subset \{P_r: \dot{X_i}=0, X_i=X_j, \Theta_i=\Theta_j,  (X_i-X_i^*)B_i=-(Z_i-Z_i^*), A_i(X_i-X_i^*)=Y_i-Y_i^*\}. $	
		 It follows from the invariant principle of the nonsmooth system \cite{shevitz1994lyapunov,bacciotti1999stability} that
		$ P_r(t)=\mathrm{col}\{\bar{X}(t), Y(t), Z(t), W(t), \Lambda(t), \Upsilon(t), \Theta(t)\}\to \mathcal{M}_r, $ as $ t\to \infty  $ and $ \mathcal{M}_r $ is the largest weakly positively invariant subset of $ \mathcal{\bar{R}}_r. $
		
		Assume that $ \hat{P}_r(t)\in \mathcal{M}_r $ for all $ t\ge 0 $ is a trajectory of \eqref{sysal:regular}. Then we have
		\begin{itemize}
			\item[a)] $ \dot{\hat{X}}_i(t)=0_{m\times r}. $
			\item[b)] $ \dot{\hat{\Lambda}}_i(t)=0_{m\times r}, $ because $ \hat{X}_i(t)=\hat{X}_j(t) $ and $ \dot{\hat{X}}_i(t)=0_{m\times r} $ .
			\item[c)] $ \dot{\hat{W}}_i(t)=0_{m\times r}, $ because $ \hat{\Theta}_i(t)=\hat{\Theta}_j(t). $
			\item[d)] Since $ \dot{\hat{X}}_i(t)=0, $ and $ A_i\hat{X}_i(t)=\hat{Y}_i(t) $ due to $ A_i(\hat{X}_i(t)-\hat{X}_i^*)=\hat{Y}_i(t)-\hat{Y}_i^* $ and $ A_i\hat{X}_i^*=\hat{Y}_i^* $, there holds $ \dot{\hat{\Upsilon}}_i(t)=0_{m_i\times r}. $
			\item[e)] Note that $ \hat{Y}_i(t)=A_i\hat{X}_i(t)=A_i\hat{X}_i(0) $ is a constant matrix. $ \dot{\hat{Y}}_i(t)=0_{m_i\times r}. $
			\item [f)]  Due to $ (\hat{X}_i(t)-\hat{X}_i^*)B_i=-(\hat{Z}_i(t)-\hat{Z}_i^*), $ and $ \dot{\hat{X}}_i(t)=0_{m\times r},  \hat{Z}_i(t)=\hat{Z}_i^*-(\hat{X}_i(0)-\hat{X}_i^*)B_i $ is a constant matrix. Because of the boundedness of the solution in (\romannumeral1), $ \dot{\hat{Z}}_i(t)=0_{m\times r_i}. $		
			\item [g)]  Recall that $ \dot{\hat{\Theta}}_i(t) $ in \eqref{sysal:least-square} is a constant matrix because $ \dot{\hat{Y}}_i(t)=0_{m_i\times r}, \dot{\hat{Z}}_i(t)=0_{m\times r_i}, \dot{\hat{W}}_i(t)=0_{m\times r_i} $. Then $ \dot{\hat{\Theta}}_i(t)=0_{m\times r}.$
		\end{itemize}
		To sum up, $ \mathcal{M}_r\subset\{P_r(t): \dot{X}_i(t)=0_{m\times r}, \dot{Y}_i(t)=0_{m_i\times r},$  $ \dot{Z}_i(t)=0_{m\times r_i}, \dot{W}_i(t)=0_{m\times r}, \dot{\Lambda}_i(t)=0_{m\times r}, \dot{\Upsilon}_i(t)=0_{m_i\times r}, \dot{\Theta}_i(t)=0_{m\times r} \} $. That is, any point in $ \mathcal{M}_r $ is an equilibrium of  \eqref{sysal:regular}, which is Lyapunov stable. In view of Lemma \ref{lemma:regular-converge} and the fact that every solution is bounded, this system is globally semi-stable. In other words, for any initial value $ P_r(0) $, with Algorithm \ref{al:regular} the solution $ P_r(t) $ converges to a Lyapunov stable equilibrium of \eqref{sysal:regular}.
		
		\noindent(\romannumeral3) Evidently, due to Lemma \ref{lemma:non-anops=an-equilibrium}, for given initial value $ P_r(0), $ every trajectory $ (\bar{X}(t), Y(t), Z(t), W(t)) $ converges to an equilibrium, which is an optimal solution of problem \eqref{nco}. Thus, $ X_i^*=X_j^* $ for all $ i,j\in \mathcal{V}. $ Furthermore, the conclusion follows.	\hfill $ \Box $
		
		\section{Discussions and Numerical Simulations}
		\label{sec:simulation}
		
		In this section, we first give a brief comparison for three cases  in Table \ref{table:1}, and then give three examples for illustration of  the validity of algorithms and theorems.
		Denote
		\begin{equation*}\label{discussion-note}
		\begin{array}{l}
		F_1=\frac{1}{2}\sum_{i=1}^{n}||X_{i}B_i-C_i+Z_i||_F^2,\\
		F_2=\frac{1}{2}\sum_{i=1}^{n}||A_iX_{i}-Y_i||_F^2,
		 \\ F_3=\frac{1}{2}\sum_{i=1}^{n}\langle\sum_{j=1}^{n}a_{i,j}(X_{i}-X_{j}), X_i\rangle_F,\\ \uppercase\expandafter{\romannumeral1}: \sum_{j=1}^{n}a_{i,j}(X_{i}-X_{j})=0_{m\times r},
		\uppercase\expandafter{\romannumeral2}: A_iX_{i}=Y_i,\\ \uppercase\expandafter{\romannumeral3}: [Y_i]_\rr-[Z_i]_\cc-\sum_{j=1}^{n}a_{i,j}(W_i-W_j)=0_{m\times r}.
		\end{array}
		\end{equation*}
		\begin{table*}[htpb]\small
			\centering
			\setlength{\abovecaptionskip}{0pt}
			\setlength{\belowcaptionskip}{0pt}
			\caption{The comparison of the three cases}
			\label{table:1}
			\begin{tabular}{|c|c|c|c|}
				\hline
				& Least Squares Case & Exact Solution Case & Regularization Case \\
				\hline
				Objective Function & $ F_1 $ & $ F_1+F_2+F_3 $ &   $F_1+\alpha\sum_{i=1}^{n}g(X_i) $  \\
				\hline
				Constraints& \uppercase\expandafter{\romannumeral1}, \uppercase\expandafter{\romannumeral2},
				\uppercase\expandafter{\romannumeral3} & \uppercase\expandafter{\romannumeral3} & \uppercase\expandafter{\romannumeral1}, \uppercase\expandafter{\romannumeral2},
				\uppercase\expandafter{\romannumeral3}  \\
				\hline
				Primal Variables& $ X_i, Y_i, Z_i, W_i $  & $ X_i, Y_i, Z_i, W_i $ & $ X_i, Y_i, Z_i, W_i $ \\
				\hline
				Dual Variables& $ \Lambda_i, \Upsilon_i, \Theta_i $ & $ \Theta_i $ & $ \Lambda_i, \Upsilon_i, \Theta_i $ \\
				\hline
				Derivative Feedback& - & - & $ \dot{X}_i $ \\
				\hline
				Exchange Status& $ X_i, W_i, \Theta_i $ & $ X_i, W_i, \Theta_i $ & $ \dot{X}_i, X_i, W_i, \Theta_i $\\
				\hline
				Convergence of $ X_i $ & $X_i\xrightarrow{\text{exponentially}} X^* $  & $ X_i\xrightarrow{\text{exponentially}} X^* $ & $ X_i\rightarrow X^* $ \\
				\hline 		
			\end{tabular}
		\end{table*}

\begin{example}\label{example:ls+ex}
	Consider a 8-node network in a label set $ \mathcal{V}=\{1,\cdots,8\}, $ whose interactions form an undirected complete graph (which has $ \mathcal{E}=\{(i,j),\forall i,j\in \mathcal{V} \} $). There is a Sylvester equation that has an exact solution:
	$ AX+XB=C, A, B, C\in\R^{8\times 8} $.
	Select a random initial value and  appropriate $ a_{i,j}>0  $, and plot the estimations for the first row of $ X_i $ over time in Fig. \ref{fig:ls01x} using Algorithm \ref{al:least-square}. As a reference, the first row of the exact solution $  X^* $ is $$ \begin{array}{r}
	X^*(1,:)=[0.5183, 0.9633,	0.2613,	0.1600,\\	0.6670,	0.2552,	0.1219,	0.4817].
	\end{array} $$ Then Fig. \ref{fig:ls01x} shows that the trajectories of the first row of $ X_i $ converge to $ X^*(1,:) $. In fact, we can also use Algorithm \ref{al:exact-solution} to solve this equation. The results show that both two algorithms derive the solutions that converge to the exact solution. Denote the error of all estimations $ E(t)=\frac{1}{n}\sum\nolimits_{i=1}^{n}||X_i-X^*||_F^2. $ Then we plot the evolution of $ log(E(t)) $ in Fig. \ref{fig:ls02log}, using Algorithms \ref{al:least-square} and \ref{al:exact-solution}, respectively. Fig. \ref{fig:ls02log} verifies the exponentially convergence of proposed distributed algorithms, which is consistent with the aforementioned theorems.
	\begin{figure}[h!]
		\centering
		\includegraphics[width=0.9\linewidth]{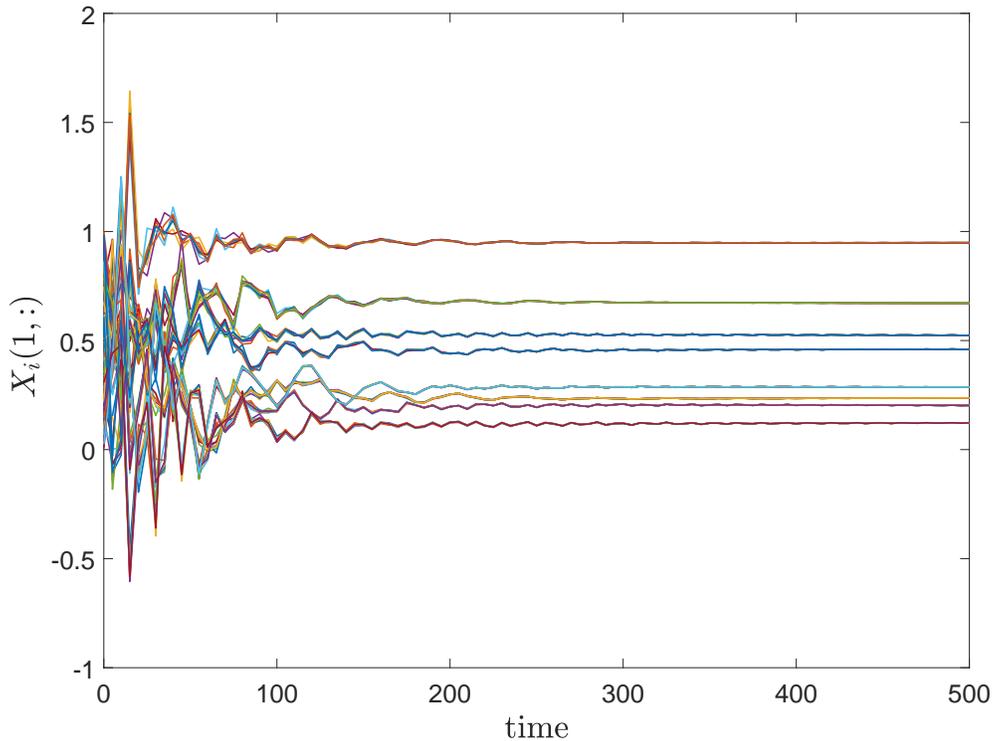}
		\caption{The estimations of all $ X_i(1,:) $ over time  using Algorithm \ref{al:least-square}.}
		\label{fig:ls01x}
	\end{figure}
\begin{figure}
	\centering
	\includegraphics[width=0.9\linewidth]{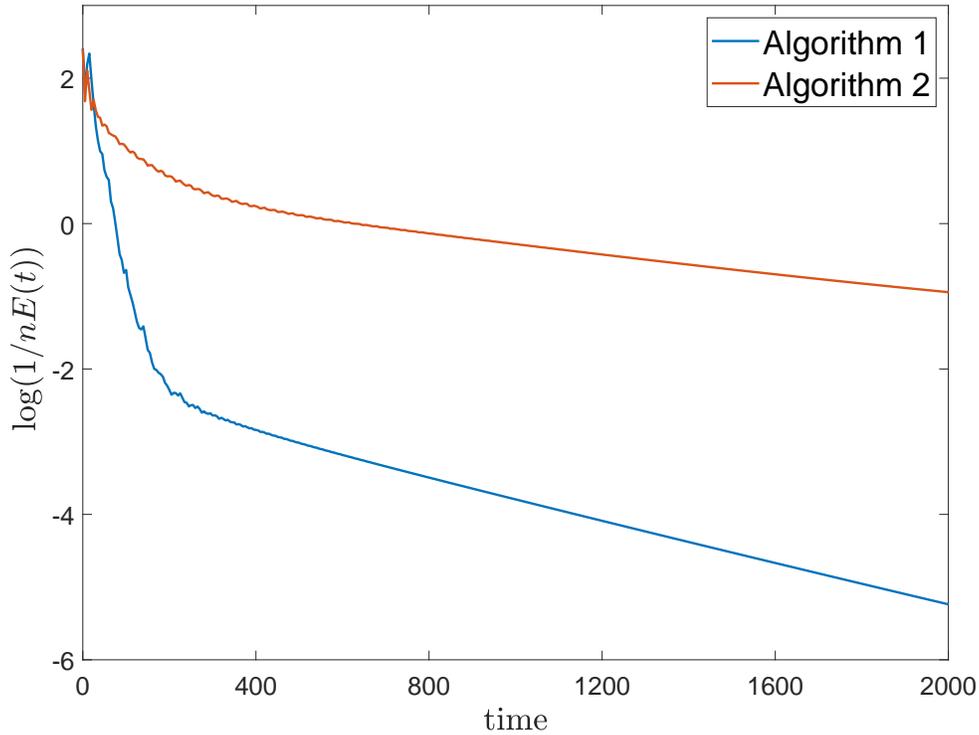}
	\caption{The evolution of $ log(E(t)) $  using Algorithms \ref{al:least-square} and \ref{al:exact-solution}, respectively.}
	\label{fig:ls02log}
\end{figure}

\end{example}	

\begin{example}\label{example:noc}
	Consider a regularization problem for the Sylvester equation
	$ AX+XB=C, A, B, C\in\R^{20\times 20} $ over a 10-node network whose interactions form an undirected complete graph. We assign two rows of $ A $ and two columns of $ B $ and $ C $ to each agent and select the regularization penalty function as $ g(X):=|X|_1=\sum_{i,j}|X_{ij}| $.
	$ l_1 $-norm of a vector is usually used to be a regularization term to get a sparse solution; similarly, we denote the sum of the absolute value of each entry of a matrix by the norm $ |X|_1. $
	Choose $ h(X) $ as a subgradient of $ g(X), h(X)\in\partial g(X), X\in \mathbb{R}^{m\times r} $ and it satisfies
	\begin{align*}
	h(X)\in \mathbb{R}^{m\times r},\quad [h(X)]_{ij}=\begin{cases}
	1,\ &if\ X_{ij}>0,\\
	0,\ &if\ X_{ij}=0,\\
	-1\ &if\ X_{ij}<0.
	\end{cases}
	\end{align*}	
	Our aim is to obtain a solution with smaller $ |\cdot|_1. $ With $ \alpha=1 $, we compare the $ |\cdot|_1 $ in Fig. \ref{fig:re01} between the exact solution and the solution trajectories by Algorithm \ref{al:regular}, and plot the error of consensus $ \frac{1}{n}\sum_{i=1}^n\sum_{j=1}^n\|X_i-X_j\|_F^2 $ in Fig. \ref{fig:re02}. Figs. \ref{fig:re01} and \ref{fig:re02} show that the solutions for all agents can reach a consensus and $ |\cdot|_1 $ of the solution has been reduced relatively.
	
	\begin{figure}
		\centering
		\includegraphics[width=0.9\linewidth]{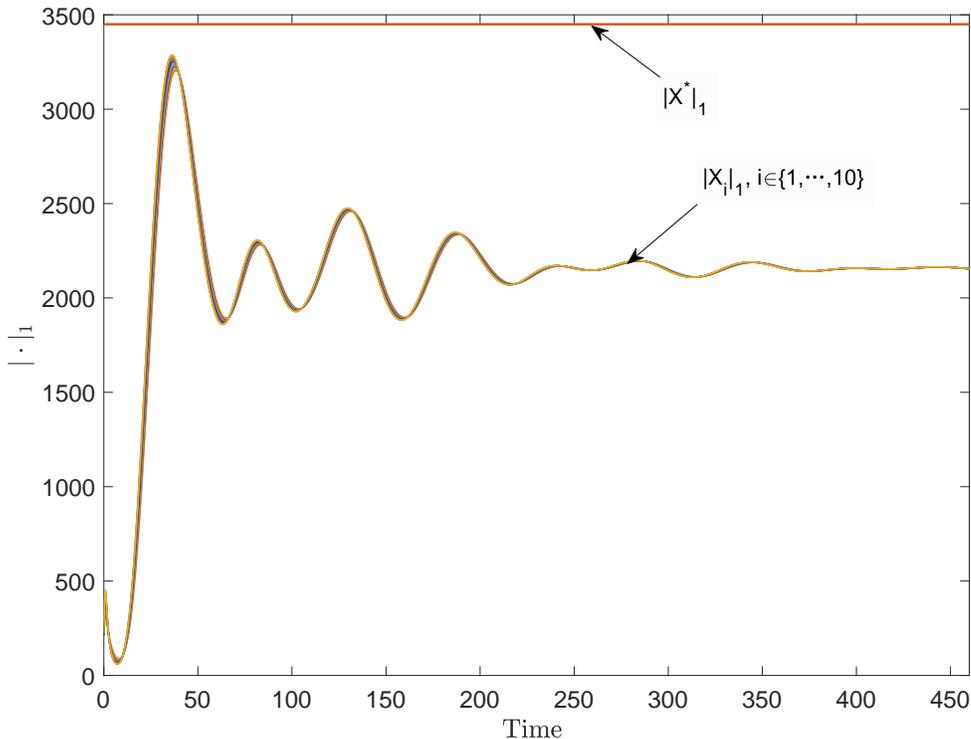}
		\caption{The evolution of $ |X^*|_1 $ and all $ |X_i(t)|_1, i\in\{1,\cdots,10\}. $}
		\label{fig:re01}
	\end{figure}
	\begin{figure}
		\centering
		\includegraphics[width=0.9\linewidth]{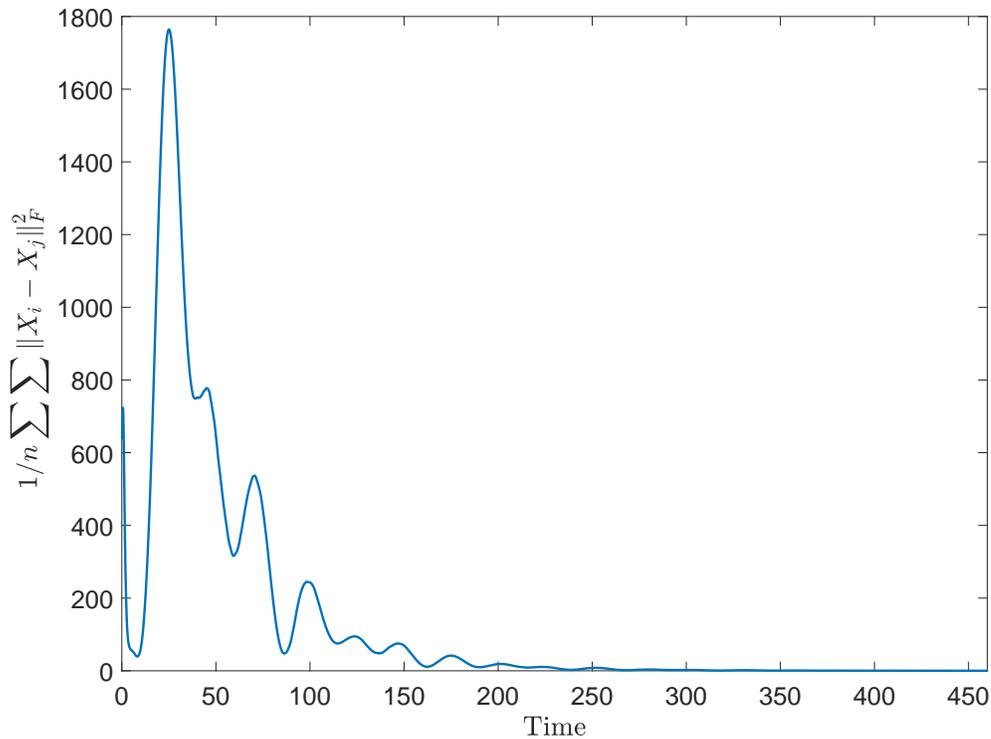}
		\caption{The error of $ 1/n\sum_{i=1}^n\sum_{j=1}^n\|X_i-X_j\|_F^2 $ over time.}
		\label{fig:re02}
	\end{figure}	
\end{example}
	\section{Conclusions}
	\label{sec:conclusion}
	This paper has proposed distributed algorithms over the multi-agent network with the help of convex optimization for solving the Sylvester equation in different cases, the least squares solution (the exact solution) case and the (nonsmooth) regularization case. In the LRRC partition  case, each agent  only has knowledge of some rows of $ A $ and some corresponding columns of $ B $ and $ C $, and holds its own state set with exchanging partial state information among its neighbors for achieving the consensus solution. Both theoretical proofs and numerical simulations have been presented to verify the convergence of proposed distributed algorithms.

\ifCLASSOPTIONcaptionsoff
  \newpage
\fi

\end{document}